\documentclass[12pt]{amsart}
\usepackage{geometry}
\usepackage{graphicx}
\usepackage{amscd}
\usepackage{amssymb}
\usepackage{amsmath}
\usepackage[latin1]{inputenc}
\usepackage{amsfonts}
\usepackage[all]{xy}

\theoremstyle{plain}                                       %
\newtheorem{thm}{\quad Theorem}                            %
\newtheorem{lem}[thm]{\quad Lemma}                         %
\newtheorem{cor}[thm]{\quad Corollary}                     %
\newtheorem{prop}[thm]{\quad Proposition}                  %
\theoremstyle{definition}                                  %
\newtheorem{defi}[thm]{\quad Definition}                   %
\newtheorem{rmk}[thm]{\quad Remark}                        %
\newtheorem{ejem}[thm]{\quad Example}                      %

\def\euler#1#2{\genfrac{\langle}{\rangle}{0pt}{}{#1}{#2}}

\newcommand{\R}{{\Bbb R}}
\newcommand{\N}{{\Bbb N}}
\newcommand{\K}{{\Bbb K}}

\newcommand{\D}{{\mathcal{D}}}

\begin{document}

\vspace{1cm}

\title{Differential equations in Ward's calculus}

\author{Ana Luz\'{o}n*, Manuel A. Mor\'{o}n$^\natural$ and Jos\'e L. Ram\'{i}rez$\dag$ }
\address{*Departamento de Matem\'{a}tica Aplicada. Universidad Polit\'{e}cnica de Madrid (Spain).}
\email{anamaria.luzon@upm.es}

\address{ $\natural$ Departamento de Algebra, Geometr\'{i}a y Topolog\'{i}a.
Universidad Complutense de Madrid  and Instituto de Matem\'{a}tica
Interdisciplinar (IMI)(Spain).} \email{mamoron@mat.ucm.es}

\address{ $\dag$ Departamento de Matem\'{a}ticas, Universidad Nacional de Colombia. Bogot\'{a}.} \email{jlramirezr@unal.edu.co}

\maketitle


\begin{abstract}
In this paper we solve some differential equations in the $D_h$ derivative in Ward's sense. We use a special metric in the formal power series ring $\K[[x]]$. The solutions of that equations are giving in terms of fixed points for certain contractive maps in our metric framework. Our main tools are Banach's Fixed Point Theorem, Fundamental Calculus Theorem and Barrow's rule for Ward's calculus. Later, we return to the usual differential calculus via Sheffer's expansion of some kind of operators. Finally, we give some examples related, in some sense, to combinatorics.
\end{abstract}

\maketitle

\section{Introduction}

Motivated by some previous works of F. H. Jackson \cite{Jack} about what nowadays is known as $q$-calculus, M. Ward introduced in \cite{ward1}, see also \cite{ward2}, an extension of the usual derivative by means of considering linear maps from polynomials into polynomials. The linear maps considered in \cite{ward1} were those that transform the polynomial $x^n$ into $h_n x^{n-1}$ for $n\in \N$, under the restrictions $h_0 =0$ and $h_n\neq0$ for $n\geq1$. Note that we recuperate the usual derivative in the special case $h_n =n$, for every $n\in \N$. Recently, see \cite{Ward}, the authors used the sequence $\{h_n\}_{n \in\N}$ in a compact way by considering the corresponding generating function $h(x)=\sum_{n\geq1}h_n x^n$. There, we associated to any such formal power series $h$ the corresponding  derivative $\D_h$. Also, we extended  to any $h$-differential calculus the pattern of relationships between Appell sequences, Sheffer sequences, and the usual derivative.  We pointed out  there, the significant role played by the $0$-Jackson derivative. The main tools, we used in \cite{Ward} to get our results, were Riordan matrices and the Hadamard product of series. As we announced in the introduction of \cite{Ward}, the current paper is focused on studying and solving differential equations in the derivative $\D_h$. To do that, in Section \ref{S:basic}, we recall our metric framework, the Banach Fixed Point Theorem, jointly with the Fundamental Calculus Theorem and Barrow's rule for Ward's calculus  \cite{Ward}. We add a brief description about the Riordan group and some of its properties that we are going to use essentially in Section \ref{S:Sheffer}.

In Section \ref{S:theory}, the existence and uniqueness of solutions for certain initial value problem (IVP) turn into fixed point problems using the Fundamental Calculus Theorem and Barrow's rule for $h$-calculus proved in \cite{Ward}. Our metric framework allows us to solve such problems by means of the Banach Fixed Point Theorem. In particular, we get an operational method, with a rigorous proof, extending to all Ward's calculus the so called Heaviside calculus for solving some types of equations. We also obtain a general method for solving the IVP for the linear  $n$-order differential equation with constant coefficients in any $h$-calculus relating this to the roots of the corresponding characteristic polynomial of the corresponding equation. To get this, we use the so called {\it reflected polynomial} as defined in \cite{Boas-Buck}. We have to say that, to get our results,  we do not use any previous results about differential equations for the classical derivative. On the contrary, our proofs give different ways to prove the classical results.

In Section \ref{S:Pascal}, we compute $\D_{h_s}$ where  $h_s(x)=x/(1-x)^s$ which are related to Pascal's triangle columns. We  obtain that, for $s\geq2$, $\D_{h_s}$ can be expressed as a finite sum in terms of the successive usual derivatives with variable coefficients. The main result needed to get it is the Vandermonde's identity.

In Section \ref{S:Sheffer}, using a Sheffer's result \cite{Shef} in the manner of Ismail \cite{Ism} in our context \cite{teo}, we write any $h$-derivative, $\D_h$,  in terms of the usual successive derivatives. We calculate the coefficients, actually monomials, going with the successive derivatives. This allows us to realize that, as occurs in the columns of Pascal's triangle, those sums that are initially infinite become finite for certain formal power series $h$. This gives us the opportunity to talk about finite and infinite differential $h$-calculus and characterize them symbolically.

In Section \ref{S:exam}, we give some examples of finite and infinite $h$-calculus with meaning in combinatorics. For these examples, we solve some concrete initial value problems.

An article related to this work has recently appeared. See \cite{VerdeStar22}. There, the authors focus on extending the Ward's derivatives to Laurent series and describe a procedure to solve some equations. Part of the significant examples in \cite{VerdeStar22}, examples 4.3, 4.5, 4.7, and 4.8 can be dealt with the tools described in our Section \ref{S:theory}.

All along this paper, we consider $\K$  is a field of characteristic zero and $\N$ is the set of natural numbers including the number zero.

\section{Previous  basic results} \label{S:basic}
\subsection{The metric framework}

Let $\K$ be a field of characteristic zero and let $\K[[x]]$ be the ring
of formal power series with coefficients in $\K$. Consider the complete ultrametric $d$ on $\K[[x]]$ given by
\[
d(f,g)=\frac{1}{2^{\omega(f-g)}}, \quad \text{for} \quad
f,\, g\in\K[[x]],
\]
where $\omega(s)$ means the \textit{order} of the power series
$s=\sum_{n=0}^{\infty}s_nx^n$, (cf.~\cite{Rob00, teo}),
defined by
\[
\omega(s)=\min\{n\in\N \ | \ s_n\neq0\}.
\]
Denote by $\cdot$ the Cauchy product of series. So,
$(\K[[x]],+,\cdot)$ has a natural structure of $\K[[x]]$-module
over the ring $(\K[[x]],+,\cdot)$. Also, $\K[[x]]$ has a structure of vector space and algebra over the field $\K$.

The metric $d$  we are using in $\K[[x]]$, all along this paper, is an ultrametric. This means that the triangular inequality  has the following strong form:
\[d(f,g)\leq \max\{d(f,s), d(s,g)\},\] for any $f,g,s \in\K[[x]]$. This is a consequence of the fact that the order in power series satisfies the property, $\omega(f+g)\geq \min\{\omega(f), \omega(g)\}$. We also know that $\omega(fg)=\omega(f)+\omega(g)$ see \cite[p. 280--281]{Rob00}.

Associated to any metric space, in particular to $(\K[[x]],d)$, we have some related concepts applicable to self-functions in $\K[[x]]$. For example, continuous functions, uniformly continuous functions, Lipschitz functions, and so on.

Consider the ultrametric space $(\K[[x]],d)$. Denote by
$End_d(\K[[x]])$ the set of all continuous endomorphisms in
$(\K[[x]],d)$ considered as a $\K$-vector space. As in the case of
classical Banach spaces, we can define what we will call the norm
associated to $d$.  We will denote it by $||\ ||_d$. We will need Corollary 26 in \cite{teo}:

\begin{cor}\label{C:d*}
The assignment $d^*:End_d(\K[[x]])\times End_d(\K[[x]])\rightarrow
\R_{+}$ given by $d^*(T_1,T_2)=||T_1-T_2||_d$ defines an ultrametric
in $End_d(\K[[x]])$.
\end{cor}

For completeness, let us recall the following definitions.
\begin{defi}
Consider the metric space $(\K[[x]],d)$ and $\Phi:(\K[[x]],d)\longrightarrow(\K[[x]],d)$ be a function, we say:
\begin{itemize}
\item[(i)] $\Phi$ is \emph{contractive}, concretely $c$-contractive, if there is a
real number $c\in[0,1)$ such that $d(\Phi(f),\Phi(g))\leq cd(f,g)$, for any $f,g\in\K[[x]]$.
\item[(ii)] $\Phi$ is \emph{non-expansive} if  $d(\Phi(f),\Phi(g))\leq d(f,g)$, for any $f,g\in\K[[x]]$.
\item[(iii)] $f\in\K[[x]]$ is a fixed point for $\Phi$ if $\Phi(f)=f$.
\end{itemize}
\end{defi}

 Note that the first two definitions are special classes of Lipschitz functions.

One of the main tool we will use along the paper is the following  well-known result.
\begin{thm}[Banach Fixed Point Theorem (BFPT)]
Let $(X,\rho)$ be a complete metric space and $f:X\rightarrow X$
contractive. Then $f$ has a unique fixed point  $x_0$ and
$f^n(x)\rightarrow x_0$ for every $x\in X$.
\end{thm}

In the above statement $f^n=f\circ f\circ\cdots\circ f$, it is the $n$-fold composition of $f$.

One of the consequences of BFPT, that we will use in the sequel, whose proof can be found in \cite[p. 1419]{convo}, is the following proposition:

\begin{prop}\label{P:f(T)}
Let $f=\sum_{n=0}^{\infty}f_nx^n$ be a power series and
$T:(\K[[x]],d)\rightarrow(\K[[x]],d)$ a contractive linear operator.
Then $f(T):(\K[[x]],d)\rightarrow(\K[[x]],d)$,
\[f(T)=\sum_{n=0}^{\infty}f_nT^n=f_0I+f_1T+f_2T^2+\cdots\] defines a
continuous  linear operator on $\K[[x]]$. Moreover
\begin{itemize}
  \item [a)] $f(T)$ is an isometry if and only if $f(0)\neq0$.
  \item [b)] $f(T)$ is contractive if and only if $f(0)=0$.
\end{itemize}
Consequently, $f(T)$ is always non-expansive.
\end{prop}

Additionally, we have the following basic properties.
\begin{prop}\label{P:propT}
Let $f,g\in\K[[x]]$ and $\lambda \in\K$. Suppose that
$T:(\K[[x]],d)\rightarrow(\K[[x]],d)$ is a contractive linear
operator.  Then
\begin{itemize}
  \item [a)] $(f+g)(T)=f(T)+g(T)$.
  \item [b)] $(\lambda f)(T)=\lambda f(T)$.
  \item [c)] $f(T)\circ g(T)=(f\cdot g)(T)$.
  \item [d)] If $g(0)=0$, then $f(g(T))=(f \circ g)(T)$.
\end{itemize}
\end{prop}
\begin{proof}
\noindent Only proofs of c) and d) are needed.  It is clear that c) is true if $f$ and $g$ are polynomials
(the same proof as in Linear Algebra of finite dimensional spaces).
To get the general result, one only has to note that the space of
polynomials, $\K[x]$, with the metric $d$ restricted to them, is dense in $(\K[[x]],d)$ and applying the definition of
the Cauchy product of series. In fact, $(\K[[x]],d)$ is the metric completion of $(\K[x],d)$.

\noindent To prove d) first note that if $g(0)=0$ and $T$ is
contractive, then $g(T)$ is also contractive. So, $f(g(T))$ makes sense, and
the composition series $f\circ g$ is well defined. And then, the
equality is clear.
\end{proof}
\subsection{The Fundamental Calculus Theorem and Barrow's rule in Ward calculus} Let us recall some definitions and basic results in Ward differential calculus that are contained in \cite{Ward}.

Suppose that $\K$  is a field of characteristic zero. Let $h(x)=\sum_{k=1}^{\infty}h_kx^k$ be a formal power series  in $\K[[x]]$ with $h_k\neq 0$  for each positive integer $k$.  The \emph{$h$-derivative matrix} $\D_h=(d_{n,m})_{n, m \in \N}$ is defined as
$$d_{n,m}=\begin{cases}
h_m, & m=n+1;\\
0, & \text{otherwise.}
\end{cases}$$
The first few rows of $\D_h$ are
\begin{align*}
\D_h=\begin{pmatrix}
 0 & h_1 & 0 & 0 & 0 \\
 0 & 0 & h_2 & 0 & 0 \\
 0 & 0 & 0 & h_3 & 0 \\
 0 & 0 & 0 & 0 & h_4 \\
 \vdots &  & \vdots & &\vdots & \ddots
 \end{pmatrix}.
 \end{align*}
 From the definition of $\D_h$ we have $\D_h(x^n)=h_nx^{n-1}$.  Therefore, if $s(x)=\sum_{k=0}^{\infty}s_kx^k$ is a formal power series, then
$$\D_h(s(x))=\sum_{k=1}^{\infty} h_ks_kx^{k-1}.$$
The \emph{$h$-integral matrix}  is defined as $\mathcal{I}_h=(i_{n,m})_{n, m \in \N}$
$$i_{n,m}=\begin{cases}
\frac{1}{h_n}, & n=m+1;\\
0, & \text{otherwise.}
\end{cases}$$
The first few rows of $\mathcal{I}_h$ are
\begin{align*}
\mathcal{I}_h=\begin{pmatrix}
 0 & 0 & 0 & 0 & 0  \\
 \frac{1}{h_1} & 0 & 0 & 0 & 0  \\
 0 & \frac{1}{h_2} & 0 & 0  & 0\\
 0 & 0 &  \frac{1}{h_3} & 0  & 0 \\
  0 & 0 & 0 &   \frac{1}{h_4} & 0 \\
 \vdots &  & \vdots & &\vdots & \ddots
 \end{pmatrix}.
 \end{align*}
 It is clear that  $\mathcal{I}_h(x^n)=\frac{1}{h_{n+1}}x^{n+1}$.  Therefore, $$\mathcal{I}_h(s(x))=\sum_{k=0}^{\infty} \frac{1}{h_{k+1}}s_{k}x^{k+1}.$$

Note that the condition $h_i\neq 0$ for all positive integer $i$ is necessary in our definition of  integral. In this paper, when we talk about $h$-differential calculus we are referring to properties related to the {\it derivative} $\D_h$  or the {\it integral} $\mathcal{I}_h$ for the power series $h$, with the needed condition $h_i\neq 0$ for all positive integer $i$. Using the products $\mathcal{I}_h\D_h$ and $\D_h\mathcal{I}_h$  we get the following general result.

\begin{thm}\label{T:TFC} The $h$-differential calculus satisfies Barrow's rule and the Fundamental Calculus Theorem for any $h$. That is, if $s=\sum_{k=0}^{\infty}s_kx^k$ is a formal power series we have
\begin{itemize}
\item Barrow's rule: $\mathcal{I}_h(\D_h(s))=s-s_0$.
\item Fundamental Theorem of $h$-Calculus: $\D_h(\mathcal{I}_h(s))=s$.
\end{itemize}
\end{thm}
The Leibniz's rule for the derivative of the product and the chain rule for the derivative of the composition are two of the main tools to compute usual derivatives of functions or formal power series. Unfortunately, in general, those rules do not hold for $h$-derivatives. In fact, an interesting characterization of the usual derivative among all $h$-derivatives, and as a counterpart of the previous general result, we have

\begin{thm} The unique derivative $\D_h$ satisfying both Leibniz's rule and chain rule is the usual derivative, i.e., $h=x/(1-x)^2.$
\end{thm}

In what follows we will use the so called {\it Hadamard product} of two formal power series, which will be denoted by $\ast$. Recall that if $f=\sum_{k=0}^{\infty}f_kx^k$ and $g=\sum_{k=0}^{\infty}g_kx^k$ are two formal power series, the Hadamard product is  defined by the formal power series $f\ast g=\sum_{k=0}^{\infty}f_{k}g_{k}x^k$. Note that the Hadamard product is associative, commutative, and distributive respect to the sum. Moreover, the power series $1/(1-x)$ is the (unique) neutral element for the Hadamard product. A power series $s=\sum_{k=0}^{\infty}s_kx^k$ is a unit for the Hadamard product if and only if $s_{k}\neq0$ for any $k\in \N$. Related to this we are going to use the following result
\begin{thm}\label{T:Had}[Theorem 12 in \cite{Ward}]
Let $s\in \K[[x]]$ and $h\in x\K[[x]]$, with $h_i\neq 0$ for all $i\geq1$. Then
$\D_h(s)=\D_0(h)*\D_0(s)=\D_0(h*s)$.
\end{thm}

\subsection{Riordan matrices and the Riordan group}
The results of this section can be found in \cite{teo} and \cite{2ways}.

\begin{defi}A \emph{Riordan matrix} is a matrix $D=(d_{i,j})_{i,j\in\N}$   whose columns are the coefficients of successive terms of a geometric progression, in $\K[[x]]$, where the initial term is a formal power series of order $0$ and the common ratio is a formal power series of order $1$.
\end{defi}
Note that for notational reasons, we recall from  \cite{teo} that a
Riordan matrix is represented as $T(f\mid
g)=D$, where $f(x)=\sum_{k=0}^{\infty}f_k x^k$ and $g(x)=\sum_{k=0}^{\infty}g_k x^k$  are formal power series in $\K[[x]]$ with $f(0)\neq0$ and $g(0)\neq0$, so that
$d_{i,j}=[x^i]x^jf(x)/g^{j+1}(x)$. Consequently, the first term is $f(x)/g(x)$ and the common ratio is  $x/g(x)$.

In this terms, Pascal's triangle is $T(1\mid 1-x)$. The above definition can be reinterpreted saying that the generating function of the $j$-th column (starting at $j=0$) of $D$ is the formal power series $x^jf(x)/g^{j+1}(x)$, which makes sense because $g(0)\neq0$. Hence, $D$ is a lower triangular matrix and it is invertible because $f(0)\neq0$.


\begin{thm} Let $T(f\mid g)=D$ be a Riordan matrix and let $\gamma(x)=\sum_{k=0}^{\infty}
\gamma_k x^k$ be a power series in $\K[[x]]$. Consider the column vector ${\bf c}=(\gamma_0, \gamma_1, \gamma_2, \dots)^T$. Then, the generating function of the matrix product $D{\bf c}$ is $\frac{f(x)}{g(x)}\gamma(\frac{x}{g(x)})$.
\end{thm}
This fact is represented by $T(f\mid g)(\gamma)=\frac{f(x)}{g(x)}\gamma(\frac{x}{g(x)})$. A proof of this result, using a special ultrametric space $(\K[[x]], d)$ can be found  in \cite[Proposition 19]{teo}.

The Riordan group (i.e., the set of all Riordan matrices) is a subgroup of the group of invertible infinite lower triangular matrices with the usual product of matrices as the
operation.

The product is given by
\[T(f\mid g)T(l\mid m)=T\left(fl\left(\frac{x}{g}\right)\big|
gm\left(\frac{x}{g}\right)\right),\] where
$fl\left(x/g\right)\equiv f(x)\cdot l\left(x/g(x)\right)$ and analogously for the second term.

The inverse is given by
\[(T(f\mid g))^{-1}\equiv T^{-1}(f\mid g)=T\left(\frac{1}{f(\frac{x}{A})}\Big| A \right),\]
where
$\left(x/A\right)\circ\left(x/g\right)=\left(x/g\right)\circ\left(x/A\right)=x$. See \cite[Proposition 20, pages 2629-2630]{teo} for more details.

The previous formal power series, denoted by $A$, is the so-called
$A$-sequence of $T(f\mid g)$. Obviously, the \emph{A-sequence} of $T(f\mid g)$ depends only on the power series $g$. Moreover, if $A=\sum_{k\geq0}a_{k}x^{k}$, then
\[d_{i,j}=\sum_{k=0}^{i-j}a_kd_{i-1,j-1+k} \qquad i,j\geq1.\]

\section{Generalities for the IVP in Ward's differential calculus}\label{S:theory}

\subsection{Some existence and uniqueness results}

Usually, $y^{(n)}$ represents the $n$-th derivative of the formal power series $y$. Consequently, from now on, we will denote by $\mathcal{D}_{h}^{(n)}$ the $n$-times composition of the operator $\mathcal{D}_h$. Consistently, we denote by $\mathcal{I}_{h}^{(n)}$ the corresponding for $\mathcal{I}_h$.

Consider the initial value problem (IVP)
\begin{equation}\label{E:IVP}
     \begin{cases}
       \mathcal{D}_{h}(y)&=G(y), \\
       y(0)&=y_0,
     \end{cases}
\end{equation}
where $G:\K[[x]]\rightarrow\K[[x]]$ is a function.

The problem (\ref{E:IVP}) could have not solutions or it could
have many of them, independently on $h$, even when $G$ is a Lipschitz map for the metric $d$. For example, considering $G(y)=\mathcal{D}_{h}(y)+p$ for $p\in\K[[x]]$
fixed. $G$ is a Lipschitz function with a Lipschitz constant equals 2, i.e., $d(G(f),G(g))\leq2d(f,g)$.
When $p\neq0$, the problem (\ref{E:IVP}) has not
any solution while if $p=0$, then any $y\in\K[[x]]$, with $y(0)=y_0$,  is
a solution.

There are many equations of the type $\mathcal{D}_{h}(y)=G(y)$. For instance, if $G(y)=q_0+q_1y+q_2y^2+\cdots+q_ny^n$ with $q_i\in\K[[x]]$, $i\in\{0,1,2,\dots,n\}$ we get what we call {\it the polynomial ordinary h-differential equation.} For $G(y)=py+q$, $p,q\in\K[[x]]$ we get {\it the linear first order h-differential equation.}

Now, we present a first existence and uniqueness result.

\begin{thm}
Let $G:(\K[[x]],d)\rightarrow(\K[[x]],d)$  be a non-expansive
function. Then (\ref{E:IVP}) has a unique solution
$y\in\K[[x]]$. This solution is the unique fixed point of $F:(\K[[x]],d)\rightarrow(\K[[x]],d)$ defined by $F(f)=y_0+(\mathcal{I}_h\circ
G)(f)$.
\end{thm}

\begin{proof}
First, we prove that  $y$ is a solution of (\ref{E:IVP}) if and only if $y$ is a fixed point of $F(f)=y_0+(\mathcal{I}_h\circ G)(f)$. If $y$ is a solution of (\ref{E:IVP}) applying $\mathcal{I}_h$ to both sides of the equation we get $\mathcal{I}_h(\mathcal{D}_h(y))=(\mathcal{I}_h\circ G)(y)$. Using (i) in Theorem \ref{T:TFC} we obtain $y=y_0 + (\mathcal{I}_h\circ G)(y)$ because $y(0)=y_0$. So, $y$ is a fixed point of $F$. On the other hand, if $y=y_0 + (\mathcal{I}_h\circ G)(y)$, then  $y(0)=y_0$. Moreover, using the linearity of $\mathcal{D}_{h}$ and (ii) in Theorem \ref{T:TFC} we obtain
\[\mathcal{D}_{h}(y)=\mathcal{D}_{h}(y_0 + (\mathcal{I}_h\circ G)(y))=\mathcal{D}_{h}(y_0)+(\mathcal{D}_{h}\circ \mathcal{I}_h)(G(y))=G(y).\]
Consequently, $y$ satisfies (\ref{E:IVP}).

Since $\mathcal{I}_h$ is $\frac{1}{2}$-contractive and $G$ is
non-expansive we get
\[
d(F(f),F(g))=\frac{1}{2^{\omega(F(f)-F(g))}}=d((\mathcal{I}_h\circ
G)(f),(\mathcal{I}_h\circ
G)(g))\leq\frac{1}{2}d(G(f),G(g))\leq\frac{1}{2}d(f,g),
\]
then, $F$ is $\frac{1}{2}$-contractive. Now, using the completeness of $d$ and  Banach's Fixed
Point Theorem we obtain that $F$ has a unique fixed point $y$,
which is the unique solution of (\ref{E:IVP}).
\end{proof}

In fact, the proof of the previous theorem can be adapted using $n$ steps to obtain the following more general result.

 \begin{thm}
\label{main:U} Given a power series $h(x)=\sum_{i\geq 1} h_ix^i$ with $h_i\in \K$ satisfying that  $h_i\neq 0$ for all positive integer $i$.
Supponse that $n$  is a positive integer and that  $y_0,y_1,\dots,y_{n-1} \in \K$. Consider
$G:(\K[[x]],d)\rightarrow(\K[[x]],d)$ satisfying
$\displaystyle{d(G(f),G(g))\leq2^{n-1}d(f,g)}$. Then the IVP
\begin{equation}\label{E:yn}
\left\{
  \begin{array}{ll}
    \mathcal{D}_{h}^{(n)}(y)&=G(y) \\
    y(0)&=y_0 \\
    \mathcal{D}_{h}(y)(0)&=y_1 \\
    &\ \vdots \\
    \mathcal{D}_{h}^{(n-1)}(y)(0)&=y_{n-1}
  \end{array}
\right.
\end{equation}
\noindent has a unique solution, which is the unique fixed point of
the contractive function
\begin{equation}\label{E:F}
    F(f)=y_0+\frac{y_1}{h_1}x+\frac{y_2}{h_{2}h_{1}}x^2+\cdots+\frac{y_{n-1}}{h_{n-1}h_{n-2}\cdots h_{1}}x^{n-1}+\mathcal{I}_{h}^{(n)}(G(f)), \quad f\in\K[[x]].
\end{equation}
 \end{thm}

\begin{rmk}
Note that for every $f\in\K[[x]]$, the sequence $\{F^n(f)\}$
converges to the  solution $y$ in $(\K[[x]],d)$.
\end{rmk}

Recall that  $y^k$ denotes the $k$-th power of $y$ respect to the Cauchy product. In the following corollary we get some significant collections of examples.

\begin{cor}
\begin{itemize}
\item[(i)] The polynomial first order h-differential equation.

Let $p_0, p_1, p_2, \dots, p_n \in \K[[x]]$. The IVP
\begin{equation}\label{E:ypn}
    \left\{
      \begin{array}{ll}
        \mathcal{D}_{h}(y)&=p_0+p_1y+p_2y^2+\cdots + p_ny^n \\
        y(0)&=y_0
      \end{array}
    \right.
\end{equation}
has a unique solution.

\item[(ii)]The linear h-differential equation of order n.

Let $q,p_0,p_1,\dots,p_{n-1}\in\K[[x]]$. The IVP
\begin{equation}\label{E:polyn}
\left\{
  \begin{array}{ll}
    \mathcal{D}_h^{(n)}(y)&=q+p_0y+p_1\mathcal{D}_{h}(y)+\cdots+p_{n-1}\mathcal{D}_{h}^{(n-1)}(y)   \\
    y(0)&=y_0\\
    \mathcal{D}_{h}(y)(0)&=y_1 \\
    & \ \vdots \\
    \mathcal{D}_{h}^{(n-1)}(y)(0)&=y_{n-1}
  \end{array}
\right.
\end{equation}
has a unique solution.

\end{itemize}
\end{cor}

\begin{proof}
It is a direct consequences of the two facts below which are straightforward to prove.

If $n$ is a positive integer and  $q,p_0,p_1,\dots, p_{n-1}, p_{n}\in\K[[x]]$, then
\begin{itemize}
\item[(i)]The function $G:(\K[[x]],d)\rightarrow(\K[[x]],d)$ defined by $G(y)=p_0+p_1y+p_2y^2+\cdots + p_ny^n $ is non-expansive.
\item[(ii)] For any  power series $h(x)=\sum_{i\geq 1} h_ix^i$ with $h_i\in \K$ satisfying that  $h_i\neq 0$ for all positive integer $i$, the function $G:(\K[[x]],d)\rightarrow(\K[[x]],d)$ defined by $G(y)=q+p_0y+p_1\mathcal{D}_{h}(y)+\cdots+p_{n-1}\mathcal{D}_{h}^{(n-1)}(y)$ is a Lipschitz function with a Lipschitz constant equal to $2^{n-1}$. \qedhere
\end{itemize}
\end{proof}

\subsection{Symbolic solutions for two particular cases}

Using the iterations given by BFPT we obtain a general expression of the solution of the IVP.

\begin{thm}\label{P:Tlineal} Consider the problem
\begin{equation}\label{E:IVPafin}
\left\{
  \begin{array}{ll}
    \mathcal{D}_{h}(y)&=T(y)+q \\
    y(0)&=y_0,
  \end{array}
\right.
\end{equation}
\noindent where $q\in\K[[x]]$ and $T:\K[[x]]\rightarrow\K[[x]]$ is
a non-expansive $\K$-endomorphism. The unique
solution $y\in\K[[x]]$  is given by
$$\frac{1}{1-x}(\mathcal{I}_{h}\circ T)(y_0+\mathcal{I}_{h}(q)).$$
\end{thm}

\begin{proof} 
The solution $y$ of (\ref{E:IVPafin}) is just the unique fixed point of the
contractive functional given by $\displaystyle{F(s)=y_0+\mathcal{I}_{h}(q)+(\mathcal{I}_{h}\circ T)(s)}$. Using now BFPT we have that for any $s\in\K[[x]]$ the sequence of iterations $\{F^n(s)\}_{n\in\N}$ converges to $y$. So, let us start
to iterate at $s_0=0$. Using the linearity of $\mathcal{I}_{h}\circ T$ one obtains
\begin{align*}
s_1&=F(0)=y_0+\mathcal{I}_{h}(q),\\
s_2&=F^2(0)=s_1+(\mathcal{I}_{h}\circ T)(s_1),\\
s_3&=F^3(0)=s_1+(\mathcal{I}_{h}\circ T)(s_1+(\mathcal{I}_{h}\circ T)(s_1))=s_1+(\mathcal{I}_{h}\circ T)(s_1)+(\mathcal{I}_{h}\circ T)^2(s_1).
\end{align*}

So, by induction on $n$
$$s_n=(I+(\mathcal{I}_{h}\circ T)+(\mathcal{I}_{h}\circ T)^2+(\mathcal{I}_{h}\circ T)^3+\cdots+(\mathcal{I}_{h}\circ T)^{n-1})(y_0+\mathcal{I}_{h}(q)).$$
 Using now the convergence in the ultrametric space $(End_d(\K[[x]]),d^*)$, see Corollary \ref{C:d*}, we obtain that \[y=\frac{I}{I-(\mathcal{I}_{h}\circ T)}(y_0+\mathcal{I}_{h}(q)). \qedhere\]
\end{proof}

Due to the use that we are going to make in the next section, we specify the following consequence of this theorem.

\begin{cor}
Consider the problem
\begin{equation}\label{E:f(V)}
\left\{
  \begin{array}{ll}
    \D_h (y)&=f(\mathcal{I}_h )(y)+ r \\
    y(0)&=y_0,
  \end{array}
\right.
\end{equation}
\noindent where $f, r\in\K[[x]]$ and $\mathcal{I}_h$ is the corresponding integral operator.
The unique solution is given by
\[
\frac{1}{1-xf}(\mathcal{I}_h)(y_0+\mathcal{I}_h (r)).
\]
\end{cor}
\begin{proof}
$f(\mathcal{I}_h )$ is a linear non-expansive operator, then from Theorem \ref{P:Tlineal} we have that the unique solution is
$$\displaystyle{y=\frac{I}{I-(\mathcal{I}_h \circ f(\mathcal{I}_h))}(y_0+\mathcal{I}_h(r))}.$$
Using now properties c) and b) in Proposition \ref{P:propT} we
obtain
\[\displaystyle{y=\frac{I}{I-xf}(\mathcal{I}_h)(y_0+\mathcal{I}_h(r))}. \qedhere\]
\end{proof}

%
%
%

\subsection{The linear $h$-differential equation of order n with constant coefficients: the Heavise
$h$-calculus.} In this section we are going to solve
 the following problem:
\begin{equation}\label{E:ynalfa}
\left\{
  \begin{array}{ll}
    \D_h ^{(n)}(y)=a_0 y +a_1\D_h (y)+\cdots+a_{n-1}\D_h ^{(n-1)}(y)+q \\
    y(0)=y_0, \ \D_h (y)(0)=y_1, \dots, \D_h ^{(n-1)}(y)(0)=y_{n-1},
  \end{array}
\right.
\end{equation}
where $a_0, a_1, \dots, a_{n-1}\in \K$ and $q\in\K[[x]]$. Moreover, we can suppose $a_0\neq0$. If it not the case, we would have a lower order problem.

We call, as usual, {\it the characteristic polynomial of the equation} in $(\ref{E:ynalfa})$ to $C(x)=x^n -\sum_{k=0}^{n-1}a_k x^k$. In these terms, the equation in problem $(\ref{E:ynalfa})$ can be restated as $C(\D_h )(y)=q$.

We apply the linear operator $\mathcal{I}_h$ to both sides in  (\ref{E:ynalfa}) and then from Barrow's rule we get
\[ \D_h ^{(n-1)}(y)=a_0 \mathcal{I}_h (y) +a_1(y - y_0 )+\cdots+a_{n-1}(\D_h ^{(n-2)}(y) - y_{n-2})+\mathcal{I}_h (q) + y_{n-1}\] or
\[ \D_h ^{(n-1)}(y)=a_0 \mathcal{I}_h (y) +a_1y  +\cdots+a_{n-1}\D_h ^{(n-2)}(y) + r_1 ,\] where $r_1 (x)=y_{n-1} -\sum_{k=1}^{n-1}a_k y_{k-1} + \mathcal{I}_h (q(x)).$ Note that $r_1 (x)$ is completely computed from data in (\ref{E:ynalfa}).
Let us write what we get if we apply only once more the integral operator to both sides in the above equality and Barrow's rule again
\[ \D_h ^{(n-2)}(y)=a_0 \mathcal{I}_h^{(2)} (y) +a_1\mathcal{I}_h (y) + a_2 y +\cdots+a_{n-1}\D_h ^{(n-3)}(y) + r_2 ,\] where $r_2 (x)=y_{n-2}-\sum_{k=2}^{n-1}a_k y_{k-2} + \mathcal{I}_h(r_1 (x))$. Again, the power series $r_2 (x)$ is completely determined from data in (\ref{E:ynalfa}). Repeating this process $(n-1)$-times, one can find an unique formal power series $r(x)$, such that the unique  solution of (\ref{E:ynalfa}) is the solution of the problem \begin{equation}\label{E: Heavi}
\left\{
  \begin{array}{ll}
  \D_h(y)&=a_0\mathcal{I}_h^{(n-1)} (y) +a_1\mathcal{I}_h^{(n-2)} (y)+\cdots+a_{n-1}y + r,   \\
    y(0)&=y_0.
  \end{array}
\right.
\end{equation}
Now, reciprocally, we can go from (\ref{E: Heavi}) to (\ref{E:ynalfa}) applying $\D_h$ $(n-1)$-times and the Fundamental Calculus Theorem in each step. With all of this we have

\begin{thm}({\it The Heaviside $h$-calculus})\label{T:h-hevy} Consider the Cauchy problem
\[\begin{cases}
C(\D_h )(y)=q,\\
y(0)=y_0, \ \D_h (y)(0)=y_1, \dots, \D_h ^{(n-1)}(y)(0)=y_{n-1}.
\end{cases}\]
Then there is a formal power series $r$, computable from data, such that the solution is given by
\[ y=\displaystyle{\frac{1}{C^* (x)}(\mathcal{I}_h)(y_0+\mathcal{I}_h(r))},\]
where $C^*(x)$ is the reflected polynomial of $C(x)$.
\end{thm}
\begin{proof}
Consider the polynomial $P(x)=\sum_{k=0}^{n-1}a_k x^k$, and $P^*$ its reflected polynomial. Note that $C(x)=x^n -P(x)$ and then $C^*(x)=1-xP^*(x)$.

As described above our problem is equivalent to $\D_h (y)=P^* (\mathcal{I}_h) +r $ with $y(0)=y_0$ for certain power series $r$. So, from the previous corollary the solution is then \[y=\displaystyle{\frac{1}{1-xP^*}(\mathcal{I}_h)(y_0+\mathcal{I}_h(r))}.\] Consequently we get the announced formula.
\end{proof}

Another important thing  we can get from Theorem \ref{T:h-hevy} is the dependence of the solutions  on the roots of the characteristic equation, $C(x)=0$, in the problem. In our conditions, $a_0\neq0$ and supposing that the field $\K$ is algebraically closed, we have that $\lambda \in \K$ is a root of $C(x)$ with multiplicity $n_\lambda$ if and only if $1/\lambda$ is a root of the reverse polynomial $C^* (x)$ with multiplicity $n_\lambda$. In fact, $C^* (x)=x^n C(1/x)$.

At this point we are going to suppose that $\K$ is algebraically
closed only for simplicity. Suppose that $\lambda_1,\dots,\lambda_s$ are
the different roots of the characteristic polynomial $C(x)$ with multiplicity
$n_1, \dots, n_s$, respectively. So, $n_1 +n_2 + \cdots + n_s =n$. Using now the decomposition into
partial  fractions we have that
\[
\frac{1}{C^* (x)}=\frac{1}{-a_0}\sum_{\ell=1}^s\left(\sum_{k=1}^{n_\ell}\frac{A_{\ell,k}}{(x-\frac{1}{\lambda_\ell})^k}\right)=
\frac{1}{a_0}\sum_{\ell=1}^s\left(\sum_{k=1}^{n_\ell}(-1)^{k+1}\lambda_\ell^k\frac{A_{\ell,k}}{(1-\lambda_\ell x)^k}\right)
.\]

So we have the following result:
\begin{cor}\label{simple}
Let $\K$ be an algebraically closed field.
Consider the Cauchy problem (\ref{E:ynalfa}) and let $C(x)=x^n -a_{n-1}x^{n-1}-\cdots -a_1 x -a_0$ be the characteristic polynomial of the equation. Suppose that $\lambda_1,\dots,\lambda_s$ are
the different roots of the characteristic polynomial $C(x)$ with multiplicity
$n_1, \dots,n_s$, respectively. Then there is a formal power series $r(x)$ computable from data such that the unique solution (\ref{E:ynalfa}) is given by
\[\frac{1}{a_0}\sum_{\ell=1}^s\left(\sum_{k=1}^{n_\ell}(-1)^{k+1}\lambda_\ell^k\frac{A_{\ell,k}}{(1-\lambda_\ell x)^k}(\mathcal{I}_h)\right)(y_0 + \mathcal{I}_h (r(x)))
.\]
\end{cor}

\section{Derivatives related to the columns of classical Pascal's Triangle.}\label{S:Pascal}

  In this section we want to point out that the derivative induced by $h_s(x)=x/(1-x)^s$, $s\geq1$,  can be written in terms of the $0$-Jackson and the usual derivatives, i.e., in terms of the derivatives induced by the first two columns of the Pascal triangle. Eventually, we prove that they can be written in terms of the successive usual derivatives but with variable coefficients.
  \begin{lem}\label{P:1/(1-x)^n}
For any positive integer $s$ and any formal power series $y$ we have
\[\D_{h_s}(y)= \sum_{k=0}^{s-1}\binom{s-1}{k}\frac{x^k}{k!}(\D_{0}(y))^{(k)},\] where $(\D_{0}(y))^{(k)}$ is the usual $k$-th derivative of the formal power series $\D_{0}(y).$
\end{lem}
\begin{proof} By Theorem \ref{T:Had} we have that $\D_{h_s}(y)= \D_{0}(h_s)\ast \D_{0}(y)$. Now \[\D_{0}(h_s)= \displaystyle{\frac{1}{(1-x)^s}}= \sum_{k\geq0}\binom{s-1+k}{k}x^k.\] For any positive integer $s$, which can be supposed to be greater than or equal to $3$, we consider the linear functional $\mathcal{H}_s:(\K[[x]], d)\longrightarrow (\K[[x]], d)$ defined by
\[\mathcal{H}_s (y)= \frac{1}{(1-x)^s}\ast y,\]
for any $y\in \K[[x]]$, where $\ast$ is the Hadamard product. Note that the matrix associated to $\mathcal{H}_s$ in the sense of \cite{teo}, as used along this paper, is the infinite diagonal matrix whose entry in the place $(m,m)$ is the binomial coefficient $\binom{s-1+m}{m}$. For every $s$, consider the linear functional
\[\mathcal{R}_s: (\K[[x]], d)\longrightarrow (\K[[x]], d),\]
defined by  $\mathcal{R}_s (y)=\sum_{k=0}^{s-1}\binom{s-1}{k}\frac{x^k}{k!}y^{(k)}$, where $y^{(k)}$ represents the usual $k$-th derivative of the power series $y$. One can easily prove that $\mathcal{R}_s$ is a non-expansive functional. To get the associated matrix we have to compute $\mathcal{R}_s(x^m )$ for every $m\geq 0$. But $$\mathcal{R}_s(x^m )=\sum_{k=0}^{s-1}\binom{s-1}{k}\binom{m}{k}x^m = \sum_{k=0}^{s-1}\binom{s-1}{k}\binom{m}{m-k}x^m.$$ Vandermonde's convolution formula implies that  $$\sum_{k=0}^{s-1}\binom{s-1}{k}\binom{m}{m-k} = \binom{s-1+m}{m}.$$ So, $\mathcal{H}_s = \mathcal{R}_s$ and consequently
\[\D_{h_s}(y)= \sum_{k=0}^{s-1}\binom{s-1}{k}\frac{x^k}{k!}(\D_{0}(y))^{(k)}.\qedhere\]
\end{proof}

\begin{thm}\label{vander} If $s\geq2$, then
\[\D_{h_s}(y)= \sum_{k=1}^{s-1} \frac{1}{k!}\binom{s-2}{k-1}x^{k-1}y^{(k)}.\]
\end{thm}

\begin{proof} Since $y=y_0 +x\D_0 (y)$, we get $y'= \D_0 (y) + x(\D_0 (y))'$. By derivating again we obtain $y''= 2(\D_0 (y))' + x(\D_0 (y))''$. We can now prove easily, by induction, that for any positive integer $k$ the equality $(\D_0 (y))^{(k)}=(y^{(k)}-k(\D_0 (y))^{(k-1)})/x$  does hold. Solving the recurrence equation we get
\[(\D_0 (y))^{(k)}=\frac{(-1)^{k}k !(y-y_0 )+\sum_{j=0}^{k-1} (-1)^j j!\binom{k}{j}x^{k-j}y^{(k-j)}}{x^{k+1}}.\] Using now Lemma \ref{P:1/(1-x)^n}, doing some computations and changing suitably the indexes of summation, we obtain
\[\D_{h_s}(y)= \sum_{k=1}^{s-1}\binom{s-1}{k}\sum_{\ell=1}^{k}(-1)^{k-\ell}\frac{1}{\ell!}x^{\ell-1}y^{(\ell)} + \D_0 (y)\sum_{k=0}^{s-1}\binom{s-1}{k}(-1)^{k}.\] The second summand above is equal to zero. Consequently
\[\D_{h_s}(y)= \sum_{k=1}^{s-1}\binom{s-1}{k}\sum_{\ell=1}^{k}(-1)^{k-\ell}\frac{1}{\ell!}x^{\ell-1}y^{(\ell)}=  \sum_{\ell=1}^{s-1}\sum_{k=\ell}^{s-1}\binom{s-1}{k}(-1)^{k-\ell}\frac{1}{\ell!}x^{\ell-1}y^{(\ell)}.\] Finally, since $$\sum_{k=\ell}^{s-1}\binom{s-1}{k}(-1)^{k-\ell}=\binom{s-2}{k-1},$$ we have the announced formula.
\end{proof}

\section{Returning to the classic differential calculus via Sheffer}\label{S:Sheffer}

In this section we extend the results of the previous one to any Ward's derivative $D_h$. Not only when $h$ is related to Pascal triangle columns. It is inspired in the work of  Mourad El-Houssieny Ismail \cite{Ism}. Sheffer's result needed is Lemma 1.1 on page 591 in \cite{Shef} restricted to the operator $D_h$:

{\it For any $h=\sum_{k=1}^{\infty}h_k x^k$ with $h_k \neq 0$ for any $k\geq 1$,  there is a sequence of polynomials $\{L^{h}_{k}\}_{k\geq 1}$,  with degree $(L^{h}_{k})\leq k-1$, such that
\begin{equation}\label{E:ShefGeneral}
\mathcal{D}_{h}(y(x))=\sum_{j\geq0}L^{h}_{j}(x)\D^{(j)}(y(x)),
\end{equation}
where $\D$ is the usual derivative.}

Note that in \cite{Ism} and \cite{Shef} expressions of the type $\sum_{k=1}^{\infty}L^{h}_{k}(x)\D^{(k)}(y(x))$ are considered as operators acting on polynomials. But it makes  sense as linear operators on $\K[[x]].$ As we know the complete ultrametric space $(\K[[x]], d)$ is the metric completion of $(\K[x], d)$ and $\D_{h}$ is always a $2$-lipschitz map on polynomials. This means that $\sum_{k=1}^{\infty}L^{h}_{k}(x)\D^{(k)}(y(x))$ is a formal power series for any formal power series $y\in \K[[x]]$.

Since $D_h(x^k)=h_kx^{k-1}$, (\ref{E:ShefGeneral}) becomes
\begin{equation}\label{E:ShefDh}
 h_k x^{k-1}=\sum_{j=1}^{k} \binom{k}{j}j!L^{h}_{j}(x)x^{k-j}.
\end{equation}
Then $L^{h}_{k}$ can be computed recurrently from (\ref{E:ShefDh}). For example, for $k=1$ it is clear that
$L_1^h(x)=h_1$; for $k=2$, we get $h_2x=\binom 21 L_1^h(x)x + \binom 22 2! L_2^h(x)=2h_1x + 2 L_2^h(x)$,
then $$L_2^h(x)=\frac x2 \left(h_2-2h_1\right).$$
For $k=3$
\begin{align*}
h_3x^2&=\binom 31 L_1^h(x)x^2 + \binom 32 2! L_2^h(x)x + \binom 33 3! L_3^h(x)\\
&=3h_1x^2 + 3\left(h_2-2h_1\right) x^2 + 6 L_3^h(x),
\end{align*}
then $$L_3^h(x)=\frac{x^2}{6} \left(3h_1-3h_2+h_3\right).$$
In general we have  the following result.
\begin{thm}\label{teo1}
Let $h(x)=\sum_{k=1}^{\infty}h_kx^k$ be in $\K[[x]]$, with $h_k\neq 0$ for all positive integer $k$.  Then
\begin{equation}\label{E:ck}
L_k^h(x)=\frac{x^{k-1}}{k!}\sum_{j=0}^k(-1)^{k+j}\binom kj h_j,
\end{equation}
for any  $k\geq0$.
\end{thm}

\begin{proof}
We will prove this using induction on $k$. From (\ref{E:ShefGeneral}), we have  $\D_h(1)=0$, then $L_0^h(x)=0$.
 We suppose that \eqref{E:ck} is true for $0\leq k\leq n$. Since $\D_h(x^{n+1})=h_{n+1}x^n$, from (\ref{E:ShefDh}) we obtain
  \[
  h_{n+1}x^n=\sum_{j=0}^{n+1} \binom{n+1}{j}j!L^{h}_{j}(x)x^{n+1-j}=(n+1)!L_{n+1}^h(x)+\sum_{j=0}^{n} \binom{n+1}{j}j!L^{h}_{j}(x)x^{n+1-j}.
  \]
By induction hypothesis
  \begin{align*}
  h_{n+1}x^n&=(n+1)!L_{n+1}^h(x)+\sum_{j=0}^{n} \binom{n+1}{j}j!x^{n+1-j}\left(\frac{x^{j-1}}{j!}\sum_{i=0}^{j}
  \binom{j}{i}(-1)^{i+j}h_i\right)\\
  &= (n+1)!L_{n+1}^h(x)+x^{n}\sum_{j=0}^{n} \sum_{i=0}^{j} \binom{n+1}{j}
  \binom{j}{i}(-1)^{i+j}h_i.
  \end{align*}
 By changing the summation indexes we obtain
  \[
  h_{n+1}x^n= (n+1)!L_{n+1}^h(x)+x^{n}\sum_{i=0}^{n} (-1)^{i} \left(\sum_{j=i}^{n} (-1)^{j} \binom{n+1}{j}
  \binom{j}{i}\right) h_i.
  \]
 Recall that  $P^{-1}P=I$ being $P$ the Pascal triangle and then $P^{-1}=\left((-1)^{k+j}\binom{k}{j}\right)_{k,j\in\N}$. Consequently, for $i\leq n$
 \[
 \sum_{j=0}^{n+1}(-1)^{n+1+j}\binom{n+1}{j}\binom{j}{i}=0,
 \]
 then
\begin{equation*}
  \sum_{j=i}^{n}(-1)^j\binom{n+1}{j}\binom{j}{i}=(-1)^n\binom{n+1}{i}.
\end{equation*}
Therefore
 \[
  h_{n+1}x^n= (n+1)!L_{n+1}^h(x)+x^{n}\sum_{i=0}^{n} (-1)^{i} (-1)^n\binom{n+1}{i}h_i,
  \]
and then
\begin{align*}
 L_{n+1}^h(x) &= \frac{x^n}{(n+1)!}\left(
 h_{n+1}+\sum_{i=0}^{n} (-1)^{n+1+i} \binom{n+1}{i}h_i\right)\\
 &=\frac{x^n}{(n+1)!}\sum_{i=0}^{n+1} (-1)^{n+1+i} \binom{n+1}{i}h_i. \qedhere
 \end{align*}
\end{proof}

From Sheffer's formula $(\ref{E:ShefGeneral})$, for any $h(x)=\sum_{n\geq1}h_nx^n$ with $h_n\neq0$  for all $n\geq1$, the operator $\D_h$ can be written in an unique way in terms of the usual derivative with variable coefficients. On the other hand, in the previous section in Theorem \ref{vander}, we described the derivatives $\D_{h_s}$, related to the columns of Pascal triangle, but as finite sums. Combining both results we infer that the polynomials $L_k^{h_s}(x)=0$,  for all $k\geq s$ in the Sheffer expression $(\ref{E:ShefGeneral})$ for $\D_{h_s}$. In particular, this implies that $\D_{h_s}(y(x))=y(x)$ becomes a linear differential equation of order $s-1$, in the usual derivatives, with variable coefficients. All of this motivates the following

\begin{defi}
  Let $h(x)=\sum_{n\geq1}h_nx^n$ with $h_n\neq0$ for all $n\geq1$. We say that the $h$-differential operator $\D_h$ is finite, or $h$ generates a finite differential calculus, if there exists $m\in \N$ such that $L_k^h(x)=0$ for all  $k>m$, where $L_k^h(x)$ are those of the Sheffer's formula $(\ref{E:ShefGeneral})$. In other case, we say that $h$ generates an infinite differential calculus.
\end{defi}

If we consider the formal power series related to the columns of Pascal's triangle, $h_s(x)=x/(1-x)^s$ then, from Theorem \ref{vander}, we have that $h_s$ \textit{generates a finite differential calculus} for $s\geq2$. In particular, $L_k^{h_s}(x)=0$ for all $k\geq s\geq 2$. Because in this case
\[\D_{h_s}(y)= \sum_{k=1}^{s-1} \frac{1}{k!}\binom{s-2}{k-1}x^{k-1}y^{(k)}.\]
Note that if $s=2$, then $\D_{h_s}=\D$.

On the other hand, if $h(x)=x/(1-x)$, then $\D_h=\D_0$. Applying Theorem \ref{teo1}, $L_k^h(x)=(-x)^{k-1}/k!$, so $h$ \textit{generates an infinite differential calculus}.

If we consider the matrix product
\[
P^{-1}h=a \Leftrightarrow
\left((-1)^{k+j}\binom{k}{j}\right)_{k,j\geq0}\left(
                                                \begin{array}{c}
                                                  h_0 \\
                                                  h_1 \\
                                                  \vdots \\
                                                  h_n \\
                                                  \vdots \\
                                                \end{array}
                                              \right)=\left(
                                                \begin{array}{c}
                                                  a_0 \\
                                                  a_1 \\
                                                  \vdots \\
                                                  a_n \\
                                                  \vdots \\
                                                \end{array}
                                              \right)
                                                \Leftrightarrow Pa=h,
\]
where $h(x)=\sum_{n\geq0}h_nx^n$, $a(x)=\sum_{n\geq0}a_nx^n$ and $P$ is the Pascal triangle or $P=T(1\mid 1-x)$ in the Riordan group notation introduced in \cite{teo}.
In particular,
\[
L_k^h(x)=\frac{x^{k-1}}{k!}\sum_{j=0}^k(-1)^{k+j}\binom kj h_j=\frac{x^{k-1}}{k!}[x^k]P^{-1}h=\frac{a_k x^{k-1}}{k!}.
\]

\begin{ejem}
  We consider $a(x)=x+x^2$, then
\[
h(x)=T(1\mid 1-x)(x+x^2)=\frac{x}{(1-x)^2}+\frac{x^2}{(1-x)^3},
\]
$h_n=[x^n]h(x)=\binom{n}{1}+\binom{n}{2}=\binom{n+1}{2}\neq0$, and
$L_0^h(x)=0$, $L_1^h(x)=1$, $L_2^h(x)=\frac{x}{2}$, $L_k^h(x)=0$ for all  $k\geq3$, then, the differential calculus generated by this $h$ is finite. In fact, it is
\[
\D_h=\D+\frac{x}{2}\D^{(2)} \ \text{or equivalently } \
\D_h(y)=y'+\frac{x}{2}y''.
\]

On the other hand, if we consider $\hat{a}(x)=x-x^2$ and making the same calculations, we can see that, in this case, the related $\hat{h}$ does not fulfill the condition $\hat{h}_n\neq0$, for all $n\geq1$, because  for $n=3$ we get $\hat{h}_3=\binom{3}{1}-\binom{3}{2}=0$.
\end{ejem}

Looking at the previous example, it is understandable to ask about some necessary and/or sufficient conditions for a polynomial $a$ makes that the related formal power series $h$, $T(1\mid 1-x)a=h$, generates a finite differential calculus. We can summarize these ideas in the following  proposition.

\begin{prop}\label{difcalcal}
  Let $h(x)=\sum_{n\geq0}h_nx^n$ and $a(x)=\sum_{n\geq0}a_nx^n$ be formal power series with $T(1\mid 1+x)h=a$. Then, $h_0=0$ and $h_n\neq0$ for all $n\geq1$ if and only if $a_0=0$ and
  \begin{equation}\label{E:BCB}
    \sum_{k=0}^{n}\binom{n}{k}a_k\neq0, \ n\geq1.
  \end{equation}
  Moreover, in this case, $h$ generates a finite differential calculus if and only if $a(x)$ is a polynomial of degree $m\geq1$. Besides,
  \begin{equation}\label{E:Dhfinito}
    \D_h=\sum_{k=1}^{m}\frac{a_k}{k!}x^{k-1}\D^{(k)}.
  \end{equation}
\end{prop}

\begin{rmk}
  Note that if in (\ref{E:BCB}) all non null coefficients are of the same sign, then $h_0=0$ and $h_n\neq0$ for all $n\geq1$. Observe that  always $a_1\neq0$.
\end{rmk}

\begin{rmk}
(\ref{E:BCB}) is related, in some sense, to the binomial coefficients bisection, shortly BCB, problem.
\end{rmk}

Collecting all the results above, we have
\begin{prop}
Let $h(x)=\sum_{n\geq0}h_nx^n$  be a formal power series. The following two facts are equivalent:
\begin{itemize}
\item[$(i)$] $h_0=0$ and $h_n\neq0$ for all $n\geq1$ and $h$ generates a finite differential calculus.
\item[$(ii)$] there exist an integer $m\geq1$ and a polynomial $a(x)=\sum_{k=1}^{m} a_{k}x^{k}$ of degree $m$, i.e., $a_{m}\neq 0$, such that $T(1\mid 1-x)a=h$
and $\sum_{k=0}^{m}\binom{n}{k}a_k\neq0$ for all $n\geq1$.
\end{itemize}
In this case  $\D_h=\sum_{k=1}^{m}\frac{a_k}{k!}x^{k-1}\D^{(k)}.$ Moreover, any  $h$ satisfying  (i), or (ii), is of the form
\[h(x)=\frac{1}{(1-x)^{m+1}}\sum_{k=1}^{m} a_{k}x^{k}(1-x)^{m-k},\]
where $m$ and $a_{k}$, for $k=1, \dots, m$, are as in (ii).
\end{prop}

\section{Some more examples of finite and infinite $h$-differential calculus}\label{S:exam}

In this section we interpret the equation  $\D_h(y)=y$  using Sheffer's expansion of $D_h$ in terms of the usual derivative. So, we are really solving initial value problems for equations in the usual derivative with variable coefficients. These problems are of finite or infinite order. We choose some examples related to combinatorics. The solutions of these problems can be given in terms of hypergeometric series. Remember that a \emph{hypergeometric series} is an expression given by
\[{_p}{F}{_q}\left(\left.\begin{tabular}{llll}$a_1,$&$a_2,$&$\dots,$&$a_p$\\$b_1,$&$b_2$,&$\dots,$&$b_q$\end{tabular}\right|x\right)=\sum_{k=0}^\infty\frac{(a_1)_k(a_2)_k\cdots(a_p)_k}{(b_1)_k(b_2)_k\cdots(b_q)_k}\frac{x^k}{k!},\]
where $(a)_n$ is the Pochamer symbol defined by $(a)_n:=a(a+1)\cdots (a+n-1)$ and $(a)_0=1$ (e.g.  see \cite{Andrew}).

\subsection{Example 1}
Consider the generating functions associated to the columns in the Pascal triangle, that is, the formal power series $h_s(x)=x/(1-x)^s$ for any positive integer $s$.
consider the initial value problem
\begin{align*}
\begin{cases}
\D_{h_s(x)}(y)&=y,\\
y(0)&=y_0.
\end{cases}
\end{align*}

For $s\geq 2$, we know that this  initial value problem is equivalent to the  equation (see   Theorem \ref{vander})
\begin{align*}
\begin{cases}
 \sum_{k=1}^{s-1} \frac{1}{k!}\binom{s-2}{k-1}x^{k-1}y^{(k)}&=y,\\
y(0)&=y_0.
\end{cases}
\end{align*}
The unique solution of this equation is the generalized exponential function (\cite{Ward})
\begin{align*}
y_0e_{h_s}^x&=y_0 \cdot {_0}{F}{_{s-2}}\left(\left.\begin{tabular}{llll}$-,$&$-,$&$\dots,$&$-$\\$2,$&$3$,&$\dots,$&$s-1$\end{tabular}\right|(s-1)!x\right) \\
&=y_0 \cdot {_1}{F}{_{s-1}}\left(\left.\begin{tabular}{llll}$ $&$ $&$1$&$ $\\$1,$&$2$,&$\dots,$&$s-1$\end{tabular}\right|(s-1)!x\right).
\end{align*}
For example, if $s=4$, the  initial value problem
$$\begin{cases}
y&=y' + xy'' + \frac{1}{6}x^2y'''\\
y(0)&=1
\end{cases}$$
has as unique solution the hypergeometric function
\begin{align*}
y&= {_1}{F}{_{3}}\left(\left.\begin{tabular}{llll}$ $&$1$&\\$1,$&$2$,&$3$\end{tabular}\right| 6x\right)\\
&=1+x+\frac{x^2}{4}+\frac{x^3}{40}+\frac{x^4}{800}+\frac{x^5}{28000}+\frac{x^6}{1568000}+\frac{x^7}{131712000}+\cdots.
   \end{align*}

\subsection{Example 2}

 In this example we consider the $h$-derivative  induced by $h_{\alpha}(x)=\sum_{k=1}^{\infty} k^{\alpha} x^k$, where $\alpha$ is a non-negative integer.  This function is related to the  polylogarithm function and for all non-negative integer $\alpha$, $h_\alpha(x)$ can be expressed as a rational function given by
 $$h_\alpha(x)=\frac{1}{(1-x)^{\alpha+1}}\sum_{i=0}^{\alpha}\euler{\alpha}{i}x^{\alpha-i},$$
where $\euler{\alpha}{k}$ are the Eulerian numbers.  For example, for $\alpha=0, 1, 2$ we have \begin{align*}
h_0(x)=\frac{x}{1-x}, \quad h_1(x)=\frac{x}{(1-x)^2}, \quad h_{2}(x)=\frac{x(x+1)}{(1-x)^3}.
\end{align*}
The action of the $h$-derivative $\D_{h_{\alpha}(x)}$  over the formal power series $s(x)=\sum_{i=0}^\infty s_ix^i$ is given by $\D_{h_{\alpha}(x)}(s(x))=\sum_{i=1}^{\infty}i^\alpha s_ix^{i-1}.$ Notice that for $\alpha=0$ and $\alpha=1$ we recover the Zero Jackson derivative and the classical derivative, respectively.

consider the initial value problem
\begin{align}\label{iniprob}
\begin{cases}
\D_{h_{\alpha}(x)}(y)&=y,\\
y(0)&=y_0.
\end{cases}
\end{align}
For this case the polynomials $L_k^h(x)$ are given by the combinatorial expression
\begin{align}\label{ecPolEx}
L_k^h(x)=\frac{x^{k-1}}{k!}\sum_{j=1}^k(-1)^{j+k}\binom kj j^\alpha, \quad \alpha\in \N.
\end{align}
If $\alpha=0$, we obtain an infinite $h$-differential operator. Indeed,
$$L_k^h(x)=\frac{x^{k-1}}{k!}\sum_{j=1}^k(-1)^{j+k}\binom kj=\frac{(-x)^{k-1}}{k!}.$$
Therefore, the initial value problem is equivalent to the following infinite equation
$$\begin{cases}
\sum_{k=1}^{\infty}\frac{(-x)^{k-1}}{k!}D^{(k)}(y)&=y,\\
y(0)&=y_0.
\end{cases}$$
The unique solution of this  system is  $y=y_0/(1-x)$. Indeed,
$$\sum_{k=1}^{\infty}\frac{(-x)^{k-1}}{k!}D^{(k)}\left(\frac{y_0}{1-x}\right)=y_0\sum_{k=1}^{\infty}\frac{(-x)^{k-1}}{k!} \frac{k!}{(1-x)^{k+1}}=y_0\sum_{k=1}^{\infty}\frac{(-x)^{k-1}}{(1-x)^{k+1}}=\frac{y_0}{1-x}.$$

If $\alpha\geq 1$, then the expression in \eqref{ecPolEx}  can be calculate by means of the Stirling number of the second kind  ${n\brace k}$. This sequence counts the number of set partitions of a set of $n$ elements into
$k$ nonempty  subsets. It is well-known the this sequence is given by the the combinatorial formula  (cf. \cite{Bernoulli})
$${n\brace k}=\frac{1}{k!}\sum_{j=0}^k\binom kj (-1)^{k+j}j^n.$$
Therefore,  $L_k^h(x)={\alpha \brace k}x^{k-1}$ and for $k>\alpha$ these polynomials are the zero polynomial. So, the initial value problem \eqref{iniprob} is equivalent to the  finite equation
$$\begin{cases}
\sum_{k=1}^{\alpha}{\alpha \brace k}x^{k-1}D^{(k)}(y)&=y,\\
y(0)&=y_0.
\end{cases}$$
The unique solution of this initial value problem is
$$y=y_0e_{h_{\alpha}(x)}^x=y_0\left(\sum_{k=0}^{\infty}\frac{x^k}{(k!)^\alpha}\right)=y_0\ {_0}{F}{_{\alpha-1}}\left(\left.\begin{tabular}{llll}$-,$&$-,$&$\dots,$&$-$\\$1,$&$1$,&$\dots,$&$1$\end{tabular}\right| x\right).$$
For example, if $\alpha=4$, we obtain the following  initial value problem
$$\begin{cases}
y&=y' + 7xy'' + 6x^2y''' + x^3y^{(4)}\\
y(0)&=y_0
\end{cases}$$
and it has as unique solution the function ($h$-exponential function)
$$y=y_0\ {_0}{F}{_{3}}\left(\left.\begin{tabular}{llll}$-,$&$-,$&$-$\\$1,$&$1$,&$1$\end{tabular}\right| x\right).$$

\subsection{Example 3: the Fibonomial calculus}
The fibonomial calculus  introduced by Krot \cite{Krot} can be obtained from the  generating function   of the Fibonacci numbers $$F(x):=\frac{x}{1-x-x^2}=\sum_{i=1}^{\infty}F_ix^i.$$
Remember that this sequence is defined by the recurrence relation $F_n=F_{n-1}+F_{n-2}$ for $n\geq 2$, with the initial conditions $F_0=0$ and $F_1=1$.

consider the initial value problem
$$\begin{cases}
\D_{F(x)}(y)&=y,\\
y(0)&=1.
\end{cases}$$
For this case the polynomials $L_k^h(x)$ are
$$L_k^h(x)=\frac{x^{k-1}}{k!}\sum_{j=1}^k(-1)^{j+k}\binom kj F_j.$$
From the combinatorial formula
$$F_n=\sum_{j=0}^n\binom{n}{j}(-1)^{j-1}F_j,$$
we get  $L_k^h(x)=\frac{(-1)^{k+1}F_k}{k!}x^{k-1}.$  Therefore the initial value problem is equivalent to the infinite functional equation:
$$\begin{cases}
\sum_{k=1}^{\infty}\frac{(-1)^{k+1}}{k!}F_k x^{k-1}D^{(k)}(y)&=y,\\
y(0)&=y_0.
\end{cases}$$
The unique solution of this equation is the Fibonacci exponential function (cf. \cite{Krot}) $e^x_{F}:=e^x_{\frac{x}{1-x-x^2}}:=\sum_{i=0}^\infty \frac{x^i}{F_i!}$, where $F_n!=F_1F_2 \cdots F_n$ and $F_0!=1$.

\subsection{Example 4: the $q$-calculus}

The $q$-derivative, also called Jackson derivative \cite{Kac}, of  formal power series is obtained by taking the function $$h_q(x)=\frac{x}{(1-x)(1-qx)}, \quad q\neq 1.$$
The $n$-th coefficient of $h_q(x)$ is given by the $n$-th $q$-integer $[n]_q=\frac{q^n-1}{q-1}=1+q+q^2+\cdots + q^{n-1}$.
consider the initial value problem
$$\begin{cases}
\D_{h_q(x)}(y)&=y,\\
y(0)&=1.
\end{cases}$$
The polynomials $L_k^h(x)$ are given by
$$L_k^h(x)=\frac{x^{k-1}}{k!}\sum_{j=1}^k(-1)^{j+k}\binom kj [j]_q.$$
From the expression
$$\sum_{j=1}^n\binom{n}{j}(-1)^{k-1}[k]_q=(q-1)^{n-1}.$$
we get the polynomials $L_k^h(x)=\frac{(q-1)^{k-1}}{k!}x^{k-1}.$  Therefore the initial value problem is equivalent to the infinite functional equation ($q\neq 1$):
$$\begin{cases}
y&=\sum_{k=1}^{\infty}\frac{(-1)^{k+1}}{k!}[k]_qx^{k-1}D^{(k)}(y)=y' + \frac{1}{2}(q-1)xy''+\frac{1}{6}(q-1)^2x^2y'''+\cdots,\\
y(0)&=1.
\end{cases}$$
The unique solution of this equation is the $q$-exponential function $e_{q}^x=\sum_{i=0}^\infty \frac{x^i}{[i]_{q}!}$.

\section{Acknowledgements}
 The first and second authors have been supported by Spanish Government Grant PID2021-126124NB-I00.

\end{document}